\documentclass[12pt]{article}
\usepackage{amssymb,amsthm,amsmath,amsfonts,latexsym,tikz,hyperref}

\newcommand{\ben}{\begin{enumerate}}
\newcommand{\een}{\end{enumerate}}
\newcommand{\ble}{\begin{lem}}
\newcommand{\ele}{\end{lem}}
\newcommand{\bth}{\begin{thm}}
\renewcommand{\eth}{\end{thm}}
\newcommand{\bpr}{\begin{prop}}
\newcommand{\epr}{\end{prop}}
\newcommand{\bco}{\begin{cor}}
\newcommand{\eco}{\end{cor}}
\newcommand{\bcon}{\begin{conj}}
\newcommand{\econ}{\end{conj}}
\newcommand{\bde}{\begin{defn}}
\newcommand{\ede}{\end{defn}}
\newcommand{\bex}{\begin{exa}}
\newcommand{\eex}{\end{exa}}
\newcommand{\barr}{\begin{array}}
\newcommand{\earr}{\end{array}}
\newcommand{\btab}{\begin{tabular}}
\newcommand{\etab}{\end{tabular}}
\newcommand{\beq}{\begin{equation}}
\newcommand{\eeq}{\end{equation}}
\newcommand{\bea}{\begin{eqnarray*}}
\newcommand{\eea}{\end{eqnarray*}}
\newcommand{\bal}{\begin{align*}}
\newcommand{\bce}{\begin{center}}
\newcommand{\ece}{\end{center}}
\newcommand{\bpi}{\begin{picture}}
\newcommand{\epi}{\end{picture}}
\newcommand{\bpp}{\begin{picture}}
\newcommand{\epp}{\end{picture}}
\newcommand{\bfi}{\begin{figure} \begin{center}}
\newcommand{\efi}{\end{center} \end{figure}}
\newcommand{\bprf}{\begin{proof}}
\newcommand{\eprf}{\end{proof}\medskip}

\newcommand{\bsl}{\begin{slide}{}}
\newcommand{\esl}{\end{slide}}
\newcommand{\bfr}{\begin{frame}}
\newcommand{\efr}{\end{frame}}

\newcommand{\hso}[1]{\hspace{-1pt}}

\newcommand{\qmq}[1]{\quad\mbox{#1}\quad}

\newcommand{\ptn}{\vdash}

\newcommand{\case}[4]{\left\{\barr{ll}#1&\mbox{#2}\\#3&\mbox{#4}\earr\right.}

\def\<{\langle}
\def\>{\rangle}

\newcommand{\ra}{\rightarrow}

\newcommand{\de}{\delta}

\newcommand{\io}{\iota}
\newcommand{\ka}{\kappa}
\newcommand{\la}{\lambda}

\newcommand{\om}{\omega}

\newcommand{\si}{\sigma}

\newcommand{\cP}{{\cal P}}

\newcommand{\cQ}{{\cal Q}}

\newcommand{\cS}{{\cal S}}

\DeclareMathOperator{\sgn}{sgn}

\newcommand{\dil}{\displaystyle}

\usepackage[hmargin=1in,vmargin=1in]{geometry}

\newtheorem{thm}{Theorem}[section]
\newtheorem{prop}[thm]{Proposition}
\newtheorem{cor}[thm]{Corollary}
\newtheorem{lem}[thm]{Lemma}
\newtheorem{conj}[thm]{Conjecture}
\newtheorem{exa}[thm]{Example}

\newcommand{\kh}{\hat{k}}
\DeclareMathOperator{\Fix}{Fix}

\begin{document}

\pagestyle{plain}

\title{Rooted partitions and number-theoretic functions
}
\author{
Bruce E. Sagan\\[-5pt]
\small Department of Mathematics, Michigan State University,\\[-5pt]
\small East Lansing, MI 48824-1027, USA, {\tt sagan@math.msu.edu}
}

\date{\today\\[10pt]
	\begin{flushleft}
	                                       \\[5pt]
	\end{flushleft}}

\maketitle

\begin{abstract}
Recently, Merca and Schmidt proved a number of identities relating partitions of an integer with two classic number-theoretic functions, namely the M\"obius function and Euler's totient function.
Their demonstrations were mainly algebraic.
We give bijective proofs of some of these results.  Our main tools are the concept of a rooted partition and an operation which we call the direct sum of a partition and a rooted partition.
\end{abstract}

\section{Introduction}

We first define three of the most well-studied functions in number theory.
Let $n$ be a nonnegative integer and
$$
[n]=\{1,2,\ldots,n\}.
$$
 A {\em partition of $n$} is a weakly decreasing sequence $\la=(\la_1,\la_2,\ldots,\la_\ell)$ of positive integers called {\em parts} satisfying $|\la|:=\sum_i\la_i =n$.
In this case we also write $\la\ptn n$.  We also let
$$
\cP(n) = \{\la \mid \la\ptn n\}
$$
so that the {\em partition function} is
$$
p(n)=\#\cP(n).
$$
Furthermore, for $n\ge 1$,  let
$$
\Phi(n) = \{k\in [n] \mid \gcd(k,n)=1\}
$$
from which we get {\em Euler's totient function}
$$
\phi(n) =\#\Phi(n)
$$
where we use the hash symbol for cardinality.
Finally, still for $n\ge1$, the {\em M\"obius function} is
\beq
\label{MuEq}
\mu(n)=\case{(-1)^{\de(n)}}{if $n$ is square free,}{0}{else,}
\eeq
where 
$$
\de(n) = \text{ number of distinct prime divisors of $n$}.
$$

In a pair of recent papers, Merca and Schmidt~\cite{MS:pir,MS:pfp} proved formulas relating $p(n)$, $\phi(n)$, and $\mu(n)$.  In order to state them, we will need some further notation.
Let
$$
S_k(n) = \text{ number of $k$'s in all the partitions of $n$}.
$$
For example, if $n=4$ and $k=1$ then
$$
\cP(4)=\{(4),\ (3,1),\ (2,2),\ (2,1,1),\ (1,1,1,1)\}
$$
and, counting the number of ones in each partition,
$$
S_1(4) = 0 + 1 + 0 + 2 + 4 = 7.
$$
These numbers had appeared earlier in the literature.
In particular, they feature in the following theorem which is due to Fine~\cite{fin:sop} although it is sometimes  referred to as Stanley's Theorem~\cite[p.\ 117, Exercise 80]{sta:ec1}.
For a thorough history of this and related results see the article of Gilbert~\cite{gil:fr}.  Note that the first sum in equation~\eqref{fg} is just $S_k(n)$.
\bth[\cite{gil:fr}]
Fix $k\ge1$  and for a partition $\la$ let 
$$
f_k(\la) = \text{ the number of parts equal to $k$ in $\la$},
$$
and
$$
g_k(\la) = \text{ the number of distinct parts of $\la$ which occur at least $k$ times}.
$$
For any $n\ge0$
\beq
\label{fg}
\sum_{\la\ptn n} f_k(\la) =\sum_{\la\ptn n} g_k(\la) .
\eeq
\eth

We will have to consider  bounding  the size of the parts of a partition.
So let
$$
S_k^{\ge r}(n) = \text{ number of $k$'s in all the partitions of $n$ with smallest part at least $r$}.
$$
Note that $S_k^{\ge1}(n) = S_k(n)$.
Also note that Merca and Schmidt used $S_{n,k}^{(r)}$ for what we are calling $S_k^{\ge r}(n)$.
But we prefer our notation as it helps to keep track of which variable is used for which parameter.
We can now state the four identities which we will prove bijectively.
The original demonstrations were mainly by manipulation of generating functions.
\bth[\cite{MS:pir,MS:pfp}]
\label{main}
For $n\ge0$ we have the following equalities.
\ben
\item $\dil S_1(n) = \sum_{k=2}^{n+1} \phi(k) S^{\ge2}_k(n+1)$.
\item $\dil p(n) = \sum_{k=3}^{n+3} \frac{\phi(k)}{2} S_k^{\ge3}(n+3)$.
\item $\dil p(n) = \sum_{k=1}^{n+1} S_k(n+1) \mu(k)$.
\item $\dil p(n) = - \sum_{k=2}^{n+2} \mu(k) S_k^{\ge2}(n+2)$.
\een
\eth

One of our crucial objects in constructing the necessary maps will be  rooted partitions.
Call a partition $\la=(\la_1,\ldots,\la_s)$ of $n$ {\em rooted} if one of its parts, say one of the $k$'s,  has been singled out.  This part is called the {\em root} and will be denoted $\kh$. 
Note that if $\la$ has $m$ parts equal to $k$ then it has $m$ corresponding partitions rooted at $k$ which are considered different and distinguished by the order of $\kh$ among the $k$'s.
For example, if $\la=(5,2,2,2,1,1)$ then there are three ways to root $\la$ at $2$, namely
$$
(5,\hat{2},2,2,1,1),\ (5,2,\hat{2},2,1,1),\ \text{and } (5,2,2,\hat{2},1,1).
$$
Let
$$
\cS_k(n) =\{\la\ptn n \mid \text{$\la$ is a rooted partition rooted at $k$}\}
$$
and, more generally
$$
\cS_k^{\ge r}(n) = \{\la\in\cS_k(n) \mid \text{$\la$ has smallest part at least $r$}\}.
$$
Clearly $\#\cS_k^{\ge r}(n) = S_k^{\ge r}(n)$ for all values of the parameters.

An important tool for combining partitions are rooted partition is what we are calling the direct sum.
Let $\la,\nu$ be two partitions with at most one of them rooted.
Their {\em direct sum} $\la\oplus\mu$ is obtained by, for each $k$, concatenating the string of $k$'s in $\la$ with the string of $k$'s in $\mu$, including the $\kh$ if one exists.
For example
$$
(5, 2, 2,  1) \oplus (4, 4, 2, \hat{2}, 2, 1, 1) = (5, 4, 4, 2, 2, 2, \hat{2}, 2, 1, 1, 1).
$$
Note that this operation is not commutative as
$$
(4, 4, 2, \hat{2}, 2, 1, 1)  \oplus (5, 2, 2,  1) = (5, 4, 4, 2, \hat{2}, 2, 2, 2, 1, 1, 1).
$$

The rest of this paper is organized as follows.
In the next section we will bijectively prove the first two equalities in Theorem~\ref{main} which relate partitions and Euler's function.
Section~\ref{cwm} will be devoted to demonstrating the last two identies between $p(n)$ and $\mu(n)$.
Most proofs will be followed by an example illustrating the maps involved.


\section{Connections with Euler's totient function}

We   begin with a bijective proof of Theorem~\ref{main} (a) which we restate here for convenience.

\bth[\cite{MS:pir}]
For $n\ge0$ we have
$$
S_1(n) = \sum_{k=2}^{n+1} \phi(k) S^{\ge2}_k(n+1).
$$
\eth
\bprf
It suffices to give a bijection $f$ from $\cS_1(n)$ to
$$
\cS'(n+1) =\{(\la',r) \mid \text{$\la'\in\cS_k^{\ge2}(n+1)$ for  some $k$ and $r\in\Phi(k)$}\}.
$$

Given $\la\in\cP_1(n)$, let $o$ be the number of ones in $\la$.  Also let $p$ be the position of $\hat{1}$, where we number the $1$ positions as $1,2,\ldots,o$ from left to right.
For ease of reading, we write 
$$
g=\gcd(o+1,p).
$$
   Define $(\la',r)=f(\la)$ as follows. 
Write
$$
\la=\mu\oplus\om.
$$
where $\om$ contains all the $1$'s in $\la$, includng $\hat{1}$, and $\mu$ is the rest of $\la$.
We now define $\om'$ to be $g$ copies of the part $(o+1)/g$ with the first of these as a root, and let
$$
\la'=\mu\oplus\om'.
$$
Finally, we let $r=p/g$.  

We must check that
$(\la',r)\in\cS'(n+1)$.
To show that $\la'\ptn n+1$, merely note that since $|\la|=n$ and $|\om|=o$ we have
$$
|\la'| =|\mu| + |\om'| = |\mu| + g\cdot (o+1)/g =|\mu|+|\om| + 1 = n+1.
$$
For the restriction on part sizes, note that $\mu$ only has parts of size at least $2$ by definition.  And since $p\le o$ it must be that $g<o+1$.  So the parts in $\om'$ are all of size $(o+1)/g\ge 2$.  Thus $\la'=\mu\oplus\om'$ has only parts greater than or equal to $2$.  Finally $r=p/g$ and the root of size $(o+1)/g$ are clearly relatively prime by definition of $g$.

To show that this map is a bijection, we construct its inverse.  Given $(\la',r)$ we write
$$
\la' = \mu \oplus \om'
$$
where $\om'$ consists of $\hat{k}$ and all the $k$'s to its right in $\la'$.  Let $s=|\om'|$ and
construct $\mu\oplus\om$ where $\om$ is $s-1$ copies of $1$.
 Finally, root the $1$ of this direct sum in position $rs/k$ to form $\la$.  It is a simple matter to check that this is a well-defined map and the inverse of $f$, so we leave those details to the reader.
\eprf

To illustrate this construction, suppose that $\la=(4,4,2,1,1,\hat{1},1,1)\in S_1(15)$.  Then we have $o=5$ ones and the position of $\hat{1}$ is $p=3$.  It follows that $g=\gcd(o+1,p)=3$.  Writing
$$
\la = (4,4,2) \oplus (1,1,\hat{1},1,1)
$$
we replace the second summand with $g=3$ copies of $(o+1)/g=2$ with the first copy as the root to get
\beq
\label{La'O}
\la'=(4,4,2) \oplus(\hat{2},2,2) = (4,4,2,\hat{2},2,2).
\eeq
For the second component of the image, we let $r=p/g=1$.

To reverse the process, suppose we are given $(\la',r)$ where $\la'$ is as in equation~\eqref{La'O} and $r=1$.
Then $\om'=(\hat{2},2,2)$ so that $k=2$ and  $s=|\om'|= 6$.  Thus $\om$ has $6-1=5$ ones and
$$
(4,4,2) \oplus (1,1,1,1,1)= (4,4,2,1,1,1,1,1).
$$
Finally we root the one in position $rs/k=3$ to obtain the original $\la$.

For the proof of Theorem~\ref{main} (b), let $n,k$ be positive integers and  consider
$$
\Psi(n) = \{ k< n/2 \mid \gcd(k,n) = 1\}.
$$
Note that $n$ and $k$ are relatively prime if and only if $n$ and $n-k$ are as well.
And for $n$ even, $n$ and $n/2$ are never coprime unless $n=2$.
It follows that for $n\ge3$ we have $\#\Psi(n) =\phi(n)/2$.
\bth[\cite{MS:pir}]
For $n\ge0$ we have
$$
p(n) = \sum_{k=3}^{n+3} \frac{\phi(k)}{2} S_k^{\ge3}(n+3).
$$
\eth
\bprf
From the discussion just before the theorem, it suffices to give a  bijection $f$ between $\cP(n)$ and
$$
\cQ^{\ge3}(n+3) =\{(\la',r) \mid \text{ $\la'\in\cS_k^{\ge3}(n+3)$ for some $k$, and $r\in\Psi(k)$}\}.
$$

Given $\la\in\cP(n)$ we define $(\la',r)=f(\la)$ as follows.  Write
$$
\la = \nu\oplus\si
$$
where $\si$ contains all the ones and twos of $\la$, and $\nu$ contains the remaining parts.
Note that if $\la$ has smallest part at least $3$ then we let $\si$ be the empty partition.
Also let 
\begin{align*}
s&= |\si|+3,\\
t&=(\text{number of twos in $\si$})+1,\\
g&=\gcd(s,t).
\end{align*}
Now construct
$$
\si'=(\kh,k,\ldots,k)
$$
where $k=s/g$ and the number of parts of $\si'$ is $g$.
Finally, define
$$
\la'=\nu\oplus\si'
$$
and 
$$
r=t/g.
$$

To verify that$f$ is well defined, we first show that $\la'\in\cQ_k^{\ge3}(n+3)$. 
The sum of the parts of $\la'$ is correct since
$$
|\la'|=|\nu|+|\si'| = |\nu| + (s/g)g = |\nu|+|\si| + 3 = n+3.
$$
We must also make sure that $\la'=\nu\oplus\si'$ has no ones or twos.  This is true of $\nu$ by definition, and for $\si'$  it suffices to show that $k\ge3$.
The number of twos in $\si$ is at most $|\si|/2$.  Now easy arithmetic manipulations imply that $s/t>2$.
It follows that $k=s/g\ge3$.

To finish the well-definedness demomstration, we must check that $r\in\Psi(k)$.
Since $r=t/g$ and $k=s/g$ it is cear that $r$ and $k$ are relatively prime.
We also need $r<k/2$.  But from the previous paragraph $s/t>2$.  This is equivalent to $kg/(rg)>2$ which is the desired inequality.

We finish by constructing $f^{-1}$.  
Suppose we are given $(\la',r)\in\cQ^{\ge3}(n+3)$
with $\la'$ rooted at $k$.  Write
$$
\la'=\nu\oplus\si'
$$
where $\si'$ contains $\kh$ and all the $k$'s after it.
Constuct a partition $\si$ of ones and twos where 
$$
|\si| =|\si'|-3
$$
and the number of twos in $\si$ is $r|\si'|/k - 1$.  
Finally, we let
$$
\la=\nu\oplus\si.
$$
We leave the demonstration that $f^{-1}$ is well defined and the inverse of $f$ to the reader.
\eprf

By way of illustration, suppose $\la=(4,4,3,2,2,1,1)$.  Then we write
$$
\la=(4,4,3,2,2,1,1) = (4,4,3) \oplus (2,2,1,1).
$$
So $\si=(2,2,1,1)$ from which we get $s=|\si|+3=9$ and $t=3$, one more than the number of twos in $\si$.
We now construct $\si'$ which is going to have $g=\gcd(s,t)=3$ copies of $k=s/g=3$ with the first one as root.
Thus
$$
\la'=(4,4,3)\oplus(\hat{3},3,3) = (4,4,3,\hat{3},3,3).
$$
Lastly, let $r=t/g=1$ to get the final pair
$$
(\la',r) = ( (4,4,3,\hat{3},3,3),\ 1).
$$

We now apply the inverse map to the pair $(\la',r)$ from the previous paragraph.  Write
$$
\la' =  (4,4,3,\hat{3},3,3) = (4,4,3)\oplus(\hat{3},3,3)
$$
so that $\si' = (\hat{3},3,3)$ and $k=3$.  We now construct $\si$ having only ones and twos with $|\si|=|\si'|-3=6$.
Since $r=1$, the number of twos  in $\si$ is to be $r|\si'|/k - 1 = 2$.  This forces $\sl=(2,2,1,1)$.  Finally
$$
\la =  (4,4,3) \oplus (2,2,1,1) = (4,4,3,2,2,1,1) 
$$
bringing us back to where we started.

\section{Connections with the M\"obius function}
\label{cwm}

In this section we will use sign-reversing involutions to construct our bijections.  We begin with a brief review of this technique.  For more information, see~\cite[Section 2.2]{sag:aoc}.

Let $\cS$ be a finite set which is {\em signed} in that it comes with  an associated function
$$
\sgn:\cS\ra\{1,-1\}.
$$
An {\em involution} on $\cS$ is a bijection $\io:\cS\ra\cS$ such that $\io^2$ is  the identity map on $\cS$.  
It is well known that an involution must decompose into cyclces $(s)$ of length one (that is, fixed points),  and $(s,t)$ of length two (that is, distinct elements $s,t$ with $\io(s)=t$ and $\io(t)=s$).
We let $\Fix\io$ be the set of fixed points of $\io$.  Call $\io$ {\em sign reversing} if
\ben
\item for each fixed point $(s)$ we have $\sgn s = 1$, and
\item for each $2$-cycle $(s,t)$ we have $\sgn s = -\sgn t$.
\een
It is clear from the definitions that
\beq
\label{io}
\sum_{s\in\cS} \sgn s = \#\Fix\io
\eeq
since all the fixed points have positive sign and the signs of elements in $2$-cycles cancel out.

As is usual when constructing sign-reversing involutions we will have two cases, depending on some parameter, which perform inverse operations which change signs.
The parameter which we will use is
$$
\pi(n)=\case{\text{smallest prime dividing $n$}}{if $n\ge2$,}{\infty}{if $n=1$,}
$$
where we consider $\infty>p$ for any prime $p$.

\bth[\cite{MS:pfp}]
\label{pS}
For $n\ge0$ we have
$$
p(n) = \sum_{k=1}^{n+1} S_k(n+1) \mu(k).
$$
\eth
\bprf
From equation~\eqref{MuEq}, we can restrict the sum in the theorem to those $k$ which are square free.  Let
$$
\cS(n) =\{ \la \mid \text{$\la\ptn n$ is a partition rooted at a square free part}\}.
$$
We turn $\cS(n)$ into a signed set by letting  the sign of a partition $\la$ with root 
$\kh$ be
$$
\sgn(\la)=\mu(k) = (-1)^{\de(k)}.
$$
Since the number of ways to root $\la$ at $k$ is the number of $k$'s in $\la$ we have 
$$
\#\cS(n+1) = \sum_k S_k(n+1)
$$
where the sum is over all square free $k$.  Futhermore, there is a bijection between the partitions $\mu$ of $n$ and the partitions of $\cS(n+1)$ obtained by inserting a $\hat{1}$ at the end of
 $\mu$. Note that the signs of the partitions in this subset are all $\mu(1)=1$ which is positive.   So, by equation~\eqref{io}, it suffices to produce a sign-reversing involution $\io$ on $\cS(n+1)$ with the rooted partitions just described as fixed points,

If $\la\in \cS(n+1)$ with root $\kh$ then let  $m$ be the number of parts equal to $k$ after and including $\kh$.
 Write
$$
\la = \mu\oplus \ka
$$
where $\ka$ contains $\kh$ and the copies of $k$ to its right.
We now have two cases for constructing $\la'=\io(\la)$ depending on the the relative sizes of $\pi(k)$ and $\pi(m)$.

Consider any $\la\in\cS(n+1)$ not ending with $\hat{1}$.  It follows that  $\min\{\pi(k),\pi(m)\}\neq\infty$ which will make the following cases well defined.  For the first case,   suppose that 
\beq
\label{Case1}
\pi(k)\le \pi(m)
\eeq 
and  let 
$$
k_1 = k/\pi(k)\qmq{and} m_1 = m \cdot \pi(k).
$$
 Now construct $\ka'$ as the partition consisting of 
 $m_1$ copies of $k_1$, where the first is rooted, and let
$$
\la' = \mu\oplus \ka'.
$$
For the second case, if we have 
$$
\pi(k)>\pi(m)
$$
then we define
$$
k_2 = k \cdot\pi(m)\qmq{and} m_2 = m/\pi(m).
$$
This time $\la'$ is constructed in exactly the same way as the first case, just using $k_2$ and $m_2$ in place of $k_1$ and $m_1$, resepctively.

We first need to prove that $\io$ is well defined in that  $\la'\in \cS(n+1)$ and $\sgn\la'=-\sgn\la$.
To show that $\la'\ptn n+1$ note that we are replacing parts summing to $km$ with parts summing to $k_1 m_1$ or $k_2 m_2$. 
But the definitions make it clear that these three  products are equal so that $\la'$ and $\la$ both partition $n+1$. We must also make sure that $k_1$ and $k_2$ are square free.
This is immediate for $k_1$ since it is a quotient of $k$.  And $k_2 = k\cdot\pi(m)$ where $\pi(m)$ is a prime smaller than the smallest prime dividing $k$.
So $k_2$ is also square free.  Finally, multiplying or dividing by a prime changes the total number of prime factors (with multiplicty) by one.
So, since we only have numbers which are square free, $\de(k)$ will also change by exactly one and thus change the sign in passing from $\la$ to $\la'$.

There remains to show that $\io^2(\la)=\la$.  Let $\la'=\io(\la)$. We have two cases resulting from the definition of $\io$.  Since the proofs are similar, we will only do the one where~\eqref{Case1} holds.
From this inequality and the fact that $k$ is square free we see that
$$
\pi(k_1) =\pi(k/\pi(k))> \pi(k) =\pi(m\cdot \pi(k)) = \pi(m_1).
$$
So $\la'$ is in the second case.  It is now trivial to verify that $\io(\la')$ yields $\la$.
\eprf

As an illustration of the previous involution, suppose $\la=(3,3,2,\hat{2},2,2,1,1)$.  Thus the root is at $k=2$ and there are $m=3$ parts of that size after and including $\hat{2}$.
Furthermore
$$
\la=(3,3,2,1,1)\oplus(\hat{2},2,2)
$$
which implies $\ka=(\hat{2},2,2)$.
Now $\pi(k)=2$ and $\pi(m)=3$ so we are in the first case of the involution.  Thus $k_1=k/\pi(k)=1$ and $m_1=m\cdot\pi(k)=6$ with corresponding partition
$\ka'=(\hat{1},1,1,1,1,1)$.  Finally
$$
\la'=(3,3,2,1,1)\oplus (\hat{1},1,1,1,1,1) =  (3,3,2,1,1,\hat{1},1,1,1,1,1).
$$

Consider what happens when we apply $\io$ to the $\la'$ just given.  Now the root is $k=1$ with $m=6$ ones after and including $\hat{1}$.
We factor
$$
\la' = (3,3,2,1,1,\hat{1},1,1,1,1,1) = (3,3,2,1,1)\oplus (\hat{1},1,1,1,1,1) 
$$
and let $\ka' =  (\hat{1},1,1,1,1,1)$.  In this situation $\pi(k)=\infty$ and $\pi(m)=2$ so we are in the second case of the involution.
Letting $k_2= k\cdot\pi(m) = 2$ and $m_2 = m/\pi(m) = 3$ yields the partition $\ka=(\hat{2},2,2)$.  Finally,
$$
\la = (3,3,2,1,1)\oplus(\hat{2},2,2) = (3,3,2,\hat{2},2,2,1,1)
$$
returning us to the original partition.

The proof of the following result is much like that of the previous theorem, so we only provide a sketch.

\bth[\cite{MS:pfp}]
For $n\ge0$ we have
$$
p(n) = - \sum_{k=2}^{n+2} \mu(k) S_k^{\ge2}(n+2).
$$
\eth
\bprf
We will prove the equality  in the form
$$
p(n) + \sum_{k=2}^{n+2} \mu(k) S_k^{\ge2}(n+2)=0
$$
So it suffices to give a sign-reversing involution $\io$ without fixed points on the set
$$
\cQ(n) = \cP(n) \uplus \{\la  \mid \text{$\la\in\cS_k^{\ge2}(n+2)$  for some $k$ which  is square free}\}
$$
where 
$$
\sgn\la = \case{1}{if $\la\in\cP(n)$,}{(-1)^{\de(k)}}{if $\la$ is rooted at $k$.}
$$

The definition of $\io$ is so similar to that of the map used in the previous theorem that we will content ourselves with defining it and leave the rest of the details to the reader.
If $\la\in\cP(n)$ then construct $\la'=\io(\la)$ as follows.  Write
$$
\la=\mu\oplus\om
$$
where $\om$ contains the ones in $\la$ or is empty if there are none.
Let 
$$
m=|\om|+2.
$$
Construct $\om'$ which consists of $m/\pi(m)$ copies of $\pi(m)$ with the first one its root.
Finally define
$$
\la'= \mu\oplus\om'.
$$
On the other hand, if $\la\in\cS_k^{\ge2}(n+2)$ then $\la'$ is constructed exactly as in the proof of Theorem~\ref{pS}.
\eprf

\section*{Acknowledgement}  We wish to than Mircea Merca for helpful discussions and references.




\nocite{*}
\bibliographystyle{alpha}

\end{document}